%% file: agt-5-2.tex
\documentclass{gtart_h}
\input agtout

\lognumber{2}
\volumenumber{5}
\volumeyear{2005}
\papernumber{2}
\pagenumbers{23}{29}
\received{15 November 2004} 
\accepted{5 January 2005}
\published{7 January 2005\nl
Corrected: 22 June 2005}

\usepackage{amsmath,amscd,amssymb}

\newtheorem{thm}{Theorem}[section]
\newtheorem*{un-no-thm}{Theorem}
\newtheorem{cor}[thm]{Corollary}     
\newtheorem{lem}[thm]{Lemma}         

\newtheorem{bigthm}{Theorem}

\theoremstyle{definition}

\theoremstyle{definition}
\theoremstyle{definition}

\theoremstyle{remark}

\newtheorem*{acks}{Acknowledgements}

\newtheorem*{rem_no}{Remark}

\newtheorem*{convent}{Conventions}

\def\smsh{\wedge}

\def\:{\colon\thinspace}
\let\Bbb\mathbb

\let\cal\mathcal

\begin{document}
\title{Poincar\'e submersions}
\covertitle{Poincar\noexpand\'e submersions}
\asciititle{Poincare submersions}

\author{John R. Klein}
\address{Department of Mathematics, Wayne State University\\Detroit, 
MI 48202, USA}
\email{klein@math.wayne.edu}
\begin{abstract} We prove
two kinds of fibering theorems for maps $X \to P$, where $X$ and $P$ are
Poincar\'e spaces.  The special case of $P = S^1$ yields
a Poincar\'e duality analogue of the fibering
theorem of Browder and Levine.
\end{abstract}
\asciiabstract{%
We prove two kinds of fibering theorems for maps X --> P, where X
and P are Poincare spaces.  The special case of P = S^1 yields a
Poincare duality analogue of the fibering theorem of Browder and
Levine.}

\primaryclass{57P10}
\secondaryclass{55R99}
\keywords{Poincar\'e duality space, fibration}
\asciikeywords{Poincare duality space, fibration}
\maketitle

\section{Introduction}
One of the early successes of surgery theory was the fibering theorem
of Browder and Levine \cite{Browder-Levine},
which gives criteria for
when a smooth map $f\:M \to S^1$ is homotopic to a submersion.
Here $M$ is assumed to be a connected
closed, smooth manifold of dimension $\ge 6$,
and we also require $f$ to induce an isomorphism of fundamental groups.
The Browder-Levine fibering theorem then says that
$f$ is homotopic to a submersion if and only if
the homotopy groups of $M$ are
finitely generated in each degree.

The purpose of the current note is to prove fibering results
in the Poincar\'e duality category.
Note that
a submersion of closed manifolds is a smooth fiber bundle with closed
manifold fibers. Replacing the closed manifolds with
finitely dominated Poincar\'e spaces and the
fiber bundle with a fibration yields the notion of
{\it Poincar\'e submersion:} this is a map between Poincar\'e
spaces whose homotopy fibers are Poincar\'e spaces.

Our first result concerns the case when the target is
acyclic (this includes the Browder-Levine situation).
Let $X$ be a connected, finitely dominated Poincar\'e duality space
of (formal) dimension $d$ and fundamental group $\pi$. Let
$$f\: X \to B\pi$$
be the classifying map for the universal cover of $X$. We will
be assuming that the classifying space $B \pi$ is a
finitely dominated Poincar\'e space of dimension $p$.

\begin{bigthm} \label{one} Let $F$ denote the
homotopy fiber of $f$. Then $F$ is a homotopy finite
Poincar\'e duality space of dimension $d-p$ if and only if
the homotopy groups of $X$ are
finitely generated in each degree.
\end{bigthm}
\bigskip

For our second result,
let $f\: X \to P$
be a
map of orientable, finitely dominated and connected Poincar\'e duality
spaces. Assume $X$ has dimension $d$ and $P$ has dimension $p$.
We will give criteria
for deciding when the homotopy fiber $F$ of $f$
satisfies Poincar\'e duality.

Let $i\: F \to X$ be the evident map.
There is an {\it umkehr homomorphism}
$$
i^{!}_*\:H_*(X) \to H_{*-p}(F)
$$
which is defined if $p\ge 3$ or if $P$ is $1$-connected
(cf.\ \S4).
The pushforward of a fundamental class $[X] \in H_d(X)$
for $X$
with respect to $i^!$ then gives a class
$$
x_f := i^!_*([X]) \in H_{d{-}p}(F)\, .
$$
This will be our candidate for a fundamental class of $F$.

\begin{bigthm} \label{two} Assume that $f$ is $2$-connected.
Then the following
are equivalent:
\begin{enumerate}

\item $H_*(F)=0$ in sufficiently large degrees.\footnote{Correction added
June, 2005: If $X$ is not $1$-connected, one also requires the hypothesis 
that the homotopy groups of $X$ are finitely generated. I am indebted to
Jonathan Hillman for pointing out that a hypothesis was missing here.
Hillman also communicated to me  the following counterexample: 
take $X$ to be the connected sum of 
$S^5 \times S^1$ with $S^3 \times S^3$ and let $f\: X \to S^1$ classify
the universal cover. Then $\pi_3(F)$ is infinitely generated.}
\item $F$ is homotopy finite.
\item $F$ is a Poincar\'e duality space.
\end{enumerate}
If in addition $X$ is $1$-connected, then
the above are equivalent to the assertion that
\begin{enumerate}
\item[\rm(4)] the homomorphism
$$
\cap x_f \:H^*(F) \to H_{d{-}p{-}*}(F)
$$
is an isomorphism in all degrees.
\end{enumerate}
\end{bigthm}

\begin{rem_no} When $P = S^p$ is a sphere, $(1) \Rightarrow (3)$
overlaps with  \cite[lemma 1.1]{Casson}. The implication
$(2) \Rightarrow (3)$ is a consequence
of \cite[theorem B]{Klein_dualizing}.

We do not {\it a priori} assume that
Poincar\'e duality spaces satisfy a finiteness condition,
so the implication $(3) \Rightarrow (2)$ is non-trivial.
\end{rem_no}

\begin{convent} A space is {\it homotopy finite} if has the homotopy
type of a finite cell complex. A space is {\it finitely dominated} if it
is the retract of a homotopy finite space.

A {\it Poincar\'e space} of formal dimension $d$
is a space $X$ for which there exists a
pair $({\cal L},[X])$ consisting of a rank one abelian system of local
coefficients ${\cal L}$ on $X$ and a (fundamental) class $[X] \in
H_d(X;{\cal L})$ such that the cap product homomorphism
$$
\cap [X] \: H^*(X;{\cal A}) \to H_{d{-}*}(X; {\cal L}\otimes {\cal A})
$$
is an isomorphism, for all local coefficient modules  ${\cal A}$ on $X$
(cf.\ \cite{Wall_pd}, \cite{Klein_pd}). If $X$ is connected,
then it is enough to establish the isomorphism when $\cal A$
is the integral group ring of the fundamental group of $X$.
When the local system ${\cal L}$ is constant, we say that
$X$ is {\it orientable}.
We do not at assume any finiteness conditions in the
definition of Poincar\'e space appearing here.
However, in the $1$-connected case,
 homotopy finiteness is actually a consequence of
Poincar\'e duality (see \ref{duality_finite} below).
\end{convent}

\begin{acks} The author is indebted to
 Mokhtar Aouina, Graeme Segal and Andrew Ranicki
for the discussions that led to  this paper.

The author was partially supported by  NSF Grant DMS-0201695.
\end{acks}

\section{Proof of Theorem \ref{one}}

We first prove the `only if' part. Assume that
$F$ is a homotopy finite Poincar\'e space.
Since $F$ is $1$-connected and homotopy finite,
we infer that its homology is finitely generated.
Apply the  mod ${\cal C}$
Hurewicz theorem (with ${\cal C}$ = the Serre class of finitely generated
abelian groups) to see that the homotopy groups of $F$ are finitely generated
\cite[corollary 9.6.16]{Spanier}.
\medskip

We now prove the `if' part.  Note that $F$ has the homotopy type of
the universal cover of $X$, so $F$ is homotopy finite dimensional
bacause $X$ is.  By the long exact homotopy sequence and the fact that
$\pi_*(X)$ is degreewise finitely generated, we infer that $\pi_*(F)$
is degreewise finitely generated.  Since $F$ is simply connected, the
mod ${\cal C}$ Hurewicz theorem shows that the homology groups of $F$
are finitely generated.  By a result of Wall \cite{Wall_finiteness},
we see that $F$ is homotopy finite.

We now know that
each space in the homotopy fiber sequence
$$
F \to X \to B\pi
$$
is finitely dominated.
It follows directly from \cite[theorem B]{Klein_dualizing}
(see also \cite{Gottlieb})
that $F$ satisfies Poincar\'e duality and has formal dimension
$d-p$. This completes the proof of Theorem \ref{one}. \qed

\section{Duality and finiteness}

A chain complex $C$ of abelian groups
is said to be {\it dualizable} if
there is chain complex $D$ and a map
$$
d\:\Bbb Z \to C\otimes D
$$
($\otimes =$ derived tensor product, and $d$ is allowed
to be degree shifting) such that, for all $P$, we get
that the induced map of complexes
$$
\hom(C,P) \to \hom(\Bbb Z,P \otimes D)
$$
(derived $\hom$) given by
$f \mapsto (f\otimes 1_D)\circ d$ induces an isomorphism on homology,
where $1_D$ denote the identity map of $D$.

A chain map $C \to D$ is said to be a {\it weak equivalence} if
it induces an isomorphism in homology. More generally $C$ and
$D$ are said to be {\it weak equivalent} if there is a finite
sequence of weak equivalences starting at $C$ and ending at $D$.
A chain complex is {\it (chain) homotopy finite} if it is weak equivalent to a finite
chain complex, i.e., a complex of
finite rank free abelian groups with finitely many non-trivial degrees.
A chain complex is
{\it finitely dominated} if is a retract up to homotopy of a finite chain
complex. It is well-known chain complex over $\Bbb Z$ is
homotopy finite if and only if it is finitely dominated
(see \cite{Wall_finiteness}).

\begin{lem} If
$C$ is {\it dualizable}, then it is homotopy finite
over $\Bbb Z$.
\end{lem}

\begin{proof} Since $\Bbb Z$ is ``compact,'' there
exists a finite chain complex $C_0$, a map $i\:C_0 \to C$
and a map $d_0\:\Bbb Z \to C_0\otimes D$ such that
$$
\begin{CD}
\Bbb Z @> d_0 >> C_0 \otimes D @>i \otimes 1 >> C \otimes D
\end{CD}
$$
is homotopic to $d$. Consider the homotopy commutative diagram
$$
\begin{CD}
\hom(C,C) @> ({-}\otimes 1_C)\circ d > \simeq > \hom(\Bbb Z,C\otimes D)\\
@A i_* AA @AA i_* A \\
\hom(C,C_0) @>\simeq > ({-}\otimes 1_C)\circ d >
\hom(\Bbb Z,C_0\otimes D)
\end{CD}
$$
The map $d_0$ lives in the lower right corner and maps to $d$ under the right
vertical map. The map $1_C$ maps to $d$ under the top horizontal map. Since
the lower horizontal map is an equivalence, we get a map
$j\: C\to C_0$ such that $i_*(j) = j \circ i$ is homotopic
to $1_C$. We conclude that
the identity map of $C$ factors up to homotopy through the finite object $C_0$.
\end{proof}

Note now if $X^d$ is a $1$-connected space which is
equipped with a chain level fundamental class $[X]$
for which Poincar\'e duality holds, then
$C(X) = $ the singular chains on $X$ is dualizable using
the maps
$$
\begin{CD}
\Bbb Z @> [X] >> C(X) @>\text{diagonal}>> C(X \times X) \simeq C(X) \otimes C(X) \, ,
\end{CD}
$$
where the first map is the homomorphism (of degree $d$) induced
by a choice of fundamental class. By the above lemma, we
infer that $C(X)$ is homotopy finite.

A result of Wall says that a $1$-connected space is homotopy finite
if and only if its chain complex is (chain) homotopy finite
(see \cite{Wall_finiteness2}).
Hence,

\begin{cor}\label{duality_finite}
 Let $X$ be  a $1$-connected space which satisfies Poincar\'e
duality. Then $X$ is also homotopy finite.
\end{cor}

\section{The umkehr homomorphism}

According to \cite[theorem 2.4]{Wall_pd},
if $\dim P \ge 3$ is a Poincar\'e duality space,
then there is a homotopy equivalence
$$
P  \,\, \simeq \,\, P_0 \cup_\alpha D^p \, ,
$$
in which $P_0$ is a CW complex of dimension $\le p{-}1$.
If $P$ is $1$-connected, then $P_0$ has the homotopy type of
a CW complex of dimension $\le p{-}2$.
If $P$ has dimension $\le 2$, then $P \simeq S^p$, and
the above decomposition is also available.

Furthermore, once an orientation for $P$ has been chosen,
the above
cell decomposition is unique up to oriented homotopy equivalence.
>From now on, we fix an identification $P := P_0 \cup D^p$, where
$\dim P_0 \le p{-}1$.

Without loss in generality, let us assume
that $f\:X \to P$ has been converted into a Hurewicz fibration.
Let $X_0 = f^{-1}(P_0)$. Then we obtain a pushout square
$$
\begin{CD}
f^{-1}(S^{p{-}1}) @>>> X_0 \\
@VVV @VVV \\
f^{-1}(D^p) @>>> X  \, .
\end{CD}
$$
Using the homotopy lifting property, we see that
the pair $(f^{-1}(D^p),f^{-1}(S^{p{-}1}))$ has the homotopy type of
the pair $(F\times D^p,F \times S^{p{-}1})$. Taking vertical cofibers
in the diagram, we get an umkehr map
$$
\begin{CD}
i^!\:X @>>> X/X_0 =  f^{-1}(D^p)/f^{-1}(S^{p{-}1}) \simeq F_+ \smsh S^p
\end{CD}
$$
The {\it umkehr homomorphism}
$$
i^!_*\:H_*(X) \to H_{*-p}(F)
$$
is the
effect of applying singular homology to $i^!$, and using
the suspension isomorphism to perform the degree shift.
 \medskip

\section{Proof of Theorem \ref{two}}\def\rk#1{#1\qua}

\rk{$(1) \Rightarrow (2)$}
By the long  exact homotopy sequence of the fibration,
we see that $\pi_*(F)$ is degreewise finitely generated. By the mod
${\cal C}$
Hurewicz theorem, we infer that $H_*(F)$ is finitely generated.
Then $F$ is homotopy finite by  \cite{Wall_finiteness}.

\rk{$(2) \Rightarrow (3)$} Follows from
\cite[theorem B]{Klein_dualizing}.

\rk{$(3) \Rightarrow (1)$} This follows from \ref{duality_finite}.

For the remainder of the
proof of the theorem, we suppose that $X$ is $1$-connected.
Then so are $F$ and $P$.

\rk{$(3) \Rightarrow (4)$} 
It will be enough to show that the class $x_f$ is a generator of
$H_{d{-}p}(F) \cong \Bbb Z$.  By definition of $x_f$, this is
equivalent to knowing that the homomorphism
$$
i^{!}_*\: H_d(X) \to H_{d-p}(F)
$$
is of degree $\pm 1$.

This can be seen as follows: the space
$X_0$ is the pullback of the fibration $f\: X\to P$
along a CW complex $P_0$ of dimension $\le p{-}2$
(this uses the fact that $P$ is $1$-connected, cf.\ \S4).
As $F$ has formal dimension $\le d{-}p$, it is straightforward to
check that $X_0$ has the homotopy type of a CW complex of dimension
$\le d{-}2$.
Using the homotopy cofiber sequence
$$
\begin{CD}
X_0 @>>> X @>i^!>> F_+\smsh S^p
\end{CD}
$$
and the fact that the homology
of $X_0$ vanishes above degree $d-2$, we see that
$i^!$ induces an isomorphism in homology in degree $d$.

\rk{$(4) \Rightarrow (3)$} Trivial.
\qed

\Addresses\recd
\end{document}
\end

%% file: agtout.tex
\def\ifplaintex{\expandafter\ifx\csname documentclass\endcsname\relax}

\def\gtp{{\mathsurround=0pt\it $\cal G\mskip-2mu$eometry \&\ 
$\cal T\!\!$opology $\cal P\!$ublications}}  

\def\recd{{\small Received:\qua\receiveddate\ifx\reviseddate\relax
\else\qquad Revised:\qua\reviseddate\fi\par}} 

\def\lognumber#1{\def\thelognumber{#1}}
\def\volumenumber#1{\def\thevolumenumber{#1}}
\def\volumeyear#1{\def\thevolumeyear{#1}}
\def\papernumber#1{\def\thepapernumber{#1}}
\def\pagenumbers#1#2{\def\startpage{#1}\def\finishpage{#2}}
\def\published#1{\def\publishdate{#1}}

\def\received#1{\def\receiveddate{#1}}

\def\accepted#1{\def\accepteddate{#1}}
\def\asciititle#1{\def\theasciititle{#1}}
\def\covertitle#1{\def\thecovertitle{#1}}

\long\def\asciiabstract#1{\long\def\theasciiabstract{#1}}
\def\asciikeywords#1{\def\theasciikeywords{#1}}

\let\\\par\let\thelognumber\relax\let\thevolumenumber\relax
\let\thepapernumber\relax\let\thevolumeyear\relax\let\startpage\relax
\let\finishpage\relax\let\publishdate\relax\let\receiveddate\relax
\let\reviseddate\relax\let\accepteddate\relax\let\theasciititle\relax
\let\thecovertitle\relax\let\theasciiauthors\relax
\let\theasciiabstract\relax\let\theasciikeywords\relax

\let\theasciiemail\relax

\ifplaintex
\font\logobig=cmssbx10 scaled 3836
\font\logomed=cmssbx10 scaled 2557
\else
\font\logobig=cmssbx10 scaled 4200
\font\logomed=cmssbx10 scaled 2800
\fi

\long\def\makeagttitle{   
\count0=\startpage
\agt\hfill      
\hbox to 45truept{\vbox to 0pt{\vglue -13truept{\logomed A\kern -.37em{\logobig 
T}\kern -.38em G}\vss}\hss}
\break
{\small Volume \thevolumenumber\ (\thevolumeyear)
\startpage--\finishpage\nl
Published: \publishdate}

\vglue .25truein

{\parskip=0pt\leftskip 0pt plus
1fil\def\\{\par\smallskip}{\Large\bf\thetitle}\par\medskip} \vglue
0.05truein

{\parskip=0pt\leftskip 0pt plus 1fil\def\\{\par}{\sc\theauthors}
\par\medskip}%
 
\vglue 0.03truein 

{\small\leftskip 25truept\rightskip 25truept{\bf Abstract}\stdspace\theabstract

{\bf AMS Classification}\stdspace\theprimaryclass
\ifx\thesecondaryclass\relax\else; \thesecondaryclass\fi\par
{\bf Keywords}\stdspace \thekeywords\par}\vglue 7truept

}   

\ifplaintex
\font\phead=cmsl9 scaled 950
\font\pnum=cmbx10 scaled 913
\font\pfoot=cmsl9 scaled 950
\headline{\vbox to 0pt{\vskip -4.5mm\line{\small\phead\ifnum
\count0=\startpage ISSN 1472-2739 (on-line) 1472-2747 (printed)
\hfill {\pnum\folio}\else\ifodd\count0\def\\{ }%
\ifx\theshorttitle\relax\thetitle\else\theshorttitle\fi\hfill{\pnum\folio}
\else\def\\{ and }{\pnum\folio}\hfill\ifx\theshortauthors\relax\theauthors
\else\theshortauthors\fi\fi\fi}\vss}}
\footline{\vbox to 0pt{\vglue 0mm\line{\small\pfoot\ifnum\count0=\startpage
\copyright\ \gtp\hfill\else
\agt, Volume \thevolumenumber\ (\thevolumeyear)\hfill\fi}\vss}}
\else
\font=cmsl9 scaled 1050
\font\lnum=cmbx10 
\font=cmsl9 scaled 1050
\makeatletter
\def\@oddhead{{\small\lhead\ifnum\count0=\startpage ISSN 1472-2739 
(on-line) 1472-2747 (printed)\hfill {\lnum\number\count0}\else\ifodd\count0
\def\\{ }\ifx\theshorttitle\relax \thetitle \else\theshorttitle\fi\hfill
{\lnum\number\count0}\else\def\\{ and }{\lnum\number\count0}
\hfill\ifx\theshortauthors\relax 
\theauthors\else\theshortauthors\fi\fi\fi}}\def\@evenhead{\@oddhead}
\def\@oddfoot{\small\lfoot\ifnum\count0=\startpage\copyright\ \gtp\hfill\else
\agt, Volume \thevolumenumber\ (\thevolumeyear)\hfill\fi}
\def\@evenfoot{\@oddfoot}
\makeatother
\fi
\let\maketitlepage\makeagttitle
\let\makeshorttitle\maketitlepage
\let\maketitle\maketitlepage

\newwrite\gtoutfile
\long\gdef\makeheadfile{  
{\def\\{, }\def\s{ }
\immediate\openout\gtoutfile head.xxx
\immediate\write\gtoutfile{Proxy-for: \ifx\theasciiauthors\relax
\theauthors\else\theasciiauthors\fi\s<\ifx\theasciiemail\relax\theemail\else\theasciiemail\fi>}
\immediate\write\gtoutfile{\noexpand\\}
\immediate\write\gtoutfile{Authors: \ifx\theasciiauthors\relax
\theauthors\else\theasciiauthors\fi}
{\def\\{ }\immediate\write\gtoutfile{Title: \ifx\theasciititle\relax
\thetitle\else\theasciititle\fi}}
\immediate\write\gtoutfile{Subj-class: GT or SG, GR etc}
\immediate\write\gtoutfile{MSC-class: \theprimaryclass\ifx\thesecondaryclass\relax\else, \thesecondaryclass\fi}
\immediate\write\gtoutfile{Journal-ref: Algebr. Geom. Topol. \thevolumenumber\s
(\thevolumeyear) \startpage-\finishpage}
\immediate\write\gtoutfile{Comments: Published by Algebraic and
Geometric Topology at}
\immediate\write\gtoutfile{\s\s\s  http://www.maths.warwick.ac.uk/agt/AGTVol\thevolumenumber/agt-\thevolumenumber-\thepapernumber.abs.html}
\immediate\write\gtoutfile{\noexpand\\}
\immediate\write\gtoutfile{}
\ifx\theasciiabstract\relax
\immediate\write\gtoutfile{\theabstract}\else
\immediate\write\gtoutfile{\theasciiabstract}\fi
\immediate\write\gtoutfile{}
\immediate\write\gtoutfile{\noexpand\\}
\immediate\write\gtoutfile{}
\immediate\closeout\gtoutfile}}  

\def\maketitlepage{\makeagttitle\makeheadfile}
\let\makeshorttitle\maketitlepage
\let\maketitle\maketitlepage